\magnification=\magstep1
 \input amstex
\NoBlackBoxes

\documentstyle{amsppt}
\topmatter

\title{INTEGRABILITY IN FINITE TERMS AND ACTIONS OF LIE GROUPS}\endtitle

\author Askold Khovanskii\endauthor

\affil{Department of Mathematics, University of Toronto, Toronto, Canada}\endaffil

\dedicatory{dedicated to Yulii Iliashenko's 75th birthday }\enddedicatory

\email{askold\@math.toronto.edu}\endemail

\thanks{The work was partially supported by the Canadian Grant No. 156833-17. }\endthanks

\abstract{ According to  Liouville's Theorem,  an idefinite integral of an elementary function is usually not an elementary function. In this notes, we discuss  that statement and a proof of this result. The differential Galois group of the extension obtained by adjoining an integral does not determine whether the integral is an elementary function or not. Nevertheless,  Liouville's Theorem can be proved using differential Galois groups. The first step  towards such a proof was suggested by Abel. This step is related to algebraic extensions and their finite Galois groups.  A significant part of this notes is dedicated  to a second step, which deals with pure transcendent extensions and their Galois groups which are connected Lie groups. The idea of the proof goes back to J.Liouville and J.Ritt. }\endabstract


\endtopmatter
\document

\head {Introduction}
\endhead
Let $K$ be a subfield of the field of meromorphic functions on a connected domain $U$ of the complex line closed under the differentiation (i.e if $f\in K$ then $f'\in K$). Such field $K$ with the operation of differentiation $f\rightarrow f'$ provides an  example of  {\it functional differential field}.  The Liouville's Theorem suggests  conditions on a function $f$ from a function differential field $K$ which are necessary and sufficient for representability of an indefinite integral of $f$  in {\it generalized elementary functions over $K$}.

In sections 1--4 we define the notions of functional differential field $K$ and  generalized elementary functions over $K$ (we follow here the exposition from [2]). A natural definition of generalized elementary functions over $K$ (see definitions 4 and 5 below) is hard to deal with.  In particular it involves a big enough list of {\it basic elementary functions} and makes use of a  non algebraic operation  of composition of  functions. An  algebraic definition (see definition 1 below) of {\it generalized elementary extension} of $K$ uses solution of the simplest differential equations  instead of composition of functions. We explain how the natural definition can be reduced to the algebraic one.

Using the algebraic definition of  generalized elementary extension one can  generalize  the Liouvvile's Theorem for an abstract differential field  $K$, whose elements are not necessarily meromorphic functions. An exposition of this result and  references to original papers can be found in [2]. This abstract algebraic result is not  directly applicable to integrals of elementary functions of one complex variable which could be multivalued, could  have singularities and so on. For its applications to elementary functions  an extra arguments (analogous to arguments we presented in sections 1--4) are needed

The differential Galois group of the extension $K\subset K(y)$ does not contain enough information to determine if the integral $y$ belongs to a generalized  elementary extension of $K$ or not. Indeed, if the integral $y$ does not belong to $K$ then the differential Galois group of $K(y)$ over $K$ is always the same: it  is isomorphic to the additive group of complex numbers. From this fact  one can conclude that the Galois theory is not sensitive enough for proving the Liouville's Theorem.

The goal of this notes is to show that  Galois theory type arguments allow to prove the Liouville's Theorem.

The first step towards such proof was suggested by Abel. This step is related to finite algebraic extensions and their finite Galois groups.

 A second step deals with a pure transcendental extension $F$ of a functional differential field $K$, obtained by adjoining $k+n$  logarithms  and exponentials,  algebraically  independent over $K$. The differential Galois group of the extension $K\subset F$  is an $k+n$-dimensional connected commutative algebraic group $G$. It has a natural represented as a group of analytic automorphisms of an analytic variety $X$. Thus $G$ acts not only on the differential field $F$  but also on other objects such as closed 1-forms on $X$. This action plays a key role in the proof.

The idea of the proof goes back to Liouville. I came up with it trying to understand  and to comment the classical  book written by J.Ritt [1]. I am grateful to Michael Singer who invite me to write comments for a new edition of this book.

\head { II. Generalized elementary functions over functional differential fields}
\endhead
In the sections 1--4 we present  definitions and general statements related  to functional and abstract differential fields and classes of their extensions including generalized elementary extensions. We follow mainly the  presentation  from the book \cite 2. In the section 5 we state Liouville's Theorem and outline its inductive proof.

\subhead{1. Abstract differential fields}
\endsubhead
A  field $F$  is said to be
{\it a differential field} if an additive map
$a\rightarrow a'$ is fixed that satisfies the Leibnitz rule $(ab)'=a'b+ab'$.
The element $a'$ is
called the  {\it derivative} of $a$. An element $y\in F$ is called {\it a constant} if $y'=0$.
All constants in $F$ form {\it the field of constants}. We add to the definition of  differential field an extra condition that {\it the field of constants is the field of complex numbers}(for our purpose it is enough to consider fields satisfying this condition).
An element $y\in F$  is said to be: {\it an exponential} of  $a$ if $y'=a'y$; {\it an exponential of integral} of  $a$ if $y'=ay$; {\it a logarithm} of  $a$ if $y'=a'/a$; {\it an integral} of $a$ if $y'=a$. In each of these cases, $y$ is defined only up to an additive or a multiplicative complex constant.

Let $K\subset F$ be a differential subfield in $F$. An element $y$ is said to be an {\it integral over $K$} if $y'=a\in K$. An {\it exponential of integral over $K$}  {\it a logarithm over $K$} and {\it an integral over $K$} are defined similarly.

Suppose that a differential field $K$ and a set $M$ lie in some differential field $F$.
{\it The adjunction} of the set $M$ to the differential field $K$ is
the minimal differential field $K\langle M\rangle$ containing both the field $K$
and the set $M$. We will refer to the transition from $K$ to $K\langle M\rangle$ as {\it adjoining} the set $M$ to the field $K$.
\definition{Definition 1}
A differential field $F$  is said to be a {\it generalized elementary extension} of a differential field  $K$
if $K\subset F$ and there exists a chain of
differential fields $K=F_0\subset \dots\subseteq F_n\supset F$
such that $F_{i+1}=F_i<y_i>$ for every $i=0$, $\dots$, $n-1$ where $y_i$ is an exponential, a logarithm, or an algebraic element over $F_i$.

An element $a\in F$ is  {\it a generalized elementary element } over $K$, $K\subset F$, if it is contained in a certain generalized elementary elementary extension of the field $K$. The following lemma is obvious.
\enddefinition
\proclaim {Lemma 1} An extension  $K\subset F$
is a {\it generalized elementary extension
if and only if there exists a chain of
differential fields $K=F_0\subseteq \dots\subseteq F_n\supset F$
such that for every $i=0$, $\dots$, $n-1$, either
$F_{i+1}$ is a finite  extension of $F_i$, or $F_{i+1}$ is a pure transcendental extension of $F_i$  obtained by
adjoining  finitely many exponentials and logarithms over $F_i$.
\endproclaim

\subhead{2. Functional differential fields and their extensions}
\endsubhead
Let $K$ be a subfield in the field $F$ of all meromorphic functions on a connected domain $U$ of the Riemann sphere $\Bbb C^1\cup \infty$ with the fixed coordinate function $x$ on $\Bbb C^1$.
Suppose that $K$ contains all complex constants and is stable under differentiation
(i.e. if $f\in K$, then $f\,'=df/dx\in K$).
Then $K$ provides an example of a {\it functional differential field.}

Let us now give a general definition.
\definition{Definition 2} Let $U,x$ be a pair consisting of a connected Riemann surface $U$ and a non constant meromorphic function $x$ on $U$. The map $f\rightarrow df/\pi^*dx$ defines the derivation  in   the field $F$ of all meromorphic
functions on  $U$ (the ratio of two meromorphic 1-forms is a well-defined meromorphic function).
 A {\it functional differential field} is any differential subfield of $F$
(containing all complex constants).
\enddefinition

The following construction helps to {\it extend} functional differential fields.
Let $K$ be a differential subfield of the field of meromorphic functions on a
connected Riemann surface $U$
equipped with a meromorphic function $x$.
Consider any connected Riemann surface $V$ together with a
nonconstant analytic map $\pi: V\to U$.
Fix the function $\pi^*x$ on $V$.
The differential field $F$ of all meromorphic functions on $V$ with the differentiation
$\varphi'= d\varphi/\pi^*dx$ contains the differential subfield $\pi^* K$ consisting of functions
of the form $\pi^* f$, where $f\in K$.
The differential field $\pi^*K$ is isomorphic to the differential field $K$,
and it lies in the differential field $F$.
For a suitable choice of the surface $V$,
an extension of the field $\pi^*K$, which is isomorphic to $K$, can be done within
the field $F$.

Suppose that we need to extend the field $K$, say, by an integral $y$ of some function $f\in K$.
This can be done in the following way.
Consider the covering of the Riemann surface $U$ by the Riemann surface $V$ of an
indefinite integral $y$ of the form $fdx$ on the surfave $U$.
By the very definition of the Riemann surface $V$,
there exists a natural projection $\pi: V\to U$, and
the function $y$ is a single-valued meromorphic function on the surface $V$.
The differential field $F$ of meromorphic functions on $V$ with the differentiation
$\varphi'= d \varphi/\pi^*dx$ contains the element $y$ as well as the field $\pi^*K$ isomorphic to $K$.
That is why the extension $\pi^*K\langle y\rangle$ is well defined as a subfield of
the differential field $F$.
We mean this particular construction of the extension whenever we talk about extensions
of functional differential fields.
The same construction allows to adjoin
a logarithm, an exponential, an integral or an exponential of integral
of any function $f$ from a functional differential field $K$ to $K$.
Similarly, for any functions $f_1,\dots,f_n\in K$, one can adjoin a
solution $y$ of an algebraic equation $y^n +f_1y^{n-1}+\dots+f_n=0$
or all the solutions
$y_1$, $\dots$, $y_n$ of this equation to $K$
(the adjunction of all the solutions $y_1$, $\dots$, $y_n$
can be implemented on the Riemann surface of the vector-function
$\bold y= y_1$, $\dots$, $y_n$).
In the same way, for any functions $f_1,\dots,f_{n+1}\in K$, one can adjoin
the $n$-dimensional $\Bbb C$-affine space of all solutions of the
linear differential equation
$y^{(n)}+f_1y^{(n-1)}+\dots +f_ny+f_{n+1}=0$ to $K$.
(Recall that a germ of any solution of this linear differential equation admits an
analytic continuation along a path on the surface $U$ not passing through the poles
of the functions $f_1$, $\dots$, $f_{n+1}$.)

Thus, {\it all above--mentioned extensions of functional differential fields can be implemented
without leaving the class of functional differential fields}.
When talking about extensions of functional differential fields, we always mean this particular
procedure.

The differential field of all complex constants and the differential field of
all rational functions of one variable can be regarded as differential
fields of functions defined on the Riemann sphere.

\subhead {3. Classes of functions and operations on multivalued functions}
\endsubhead
An indefinite integral of an elementary function is a function rather than an element of an abstract differential field.
In functional spaces, for example, apart from differentiation and algebraic operations,
an absolutely non-algebraic operation is defined, namely, the composition.
Anyway, functional spaces provide more means for writing ``explicit formulas''
than abstract differential fields.
Besides, we should take into account that functions can be multivalued,
can have singularities and so on.

In functional spaces, it is not hard to formalize the problem of unsolvability
of equations in explicit form.
One can proceed as follows:
fix a class of functions and say that an equation is solvable explicitly if
its solution belongs to this class.
Different classes of functions correspond to different notions of solvability.

\subhead{3.1. Defining classes of functions by the lists
of data}
\endsubhead
A class of functions can be introduced by specifying a list of {\it basic functions}
and a list of {\it admissible operations}.
Given the two lists, the class of functions is defined as the set
of all functions that can be obtained
from the basic functions by repeated application of admissible operations.
Below, we define the class of
{\it generalized elementary functions} and the class of {\it generalized elementary functions over a functional differential field $K$} in exactly this way.

Classes of functions, which appear in the problems of integrability in finite terms,
contain multivalued functions.
Thus the basic terminology should be made clear.
We work with multivalued functions ``globally'', which leads
to a more general understanding of classes of functions defined by
lists of basic functions and of admissible operations.
A multivalued function is regarded as a single entity.
{\it Operations on multivalued functions} can be defined.
The result of such an operation is a set of multivalued functions;
every element of this set is called a function obtained from the given functions
by the given operation.
A  {\it class of functions} is defined as the set of all (multivalued) functions
that can be obtained from the basic functions by repeated
application of admissible operations.

\subhead{3.2. Operations on multivalued functions}
\endsubhead
Let us define, for example, the sum of two multivalued functions on a connected Riemann surface $U$.

\definition{Definition 3}
Take an arbitrary point $a$ in $U$, any germ
$f_a$ of an analytic function $f$ at the point $a$ and any germ $g_a$ of an analytic function $g$
at the same point $a$.
We say that the multivalued function $\varphi$ on $U$ generated by the germ $\varphi_a=f_a+g_a$
{\it is representable as the sum of the functions} $f$ and $g$.
\enddefinition

For example, it is easy to see that exactly two functions of one variable
are representable in the form
$\sqrt{x}+\sqrt{x}$, namely, $f_1=2\sqrt{x}$ and $f_2\equiv 0$.
Other operations on multivalued functions are defined in exactly the same way.
{\it For a class of multivalued functions, being stable under addition means that,
together with any pair
of its functions, this class contains all functions representable as their sum.}
The same applies to all other operations on multivalued functions understood
in the same sense as above.

In the definition given above, not only the operation of addition plays a key role but
also the operation of analytic continuation hidden in the notion
of multivalued function.
Indeed, consider the following example.
Let $f_1$ be an analytic function defined on an open subset $V$ of the complex line
$\Bbb C^1$ and admitting no analytic continuation outside of $V$,
and let $f_2$ be an analytic function on
$V$ given by the formula $f_2=-f_1$.
According to our definition, the zero function is representable in the form
$f_1+f_2$ {\it on the entire complex line}.
By the commonly accepted viewpoint, the equality $f_1+f_2=0$ holds inside
the region $V$ but not outside.

Working with multivalued functions globally, we do not insist on the existence of
{\it a common region}, were all necessary operations would be performed on
single-valued branches of multivalued functions.
A first operation can be performed in a first region,
then a second operation can be performed in a second, different region
on analytic continuations of functions obtained on the first step.
In essence, this more general understanding of operations is equivalent to including
analytic continuation to the list of admissible operations on the analytic germs.
\subhead{4. Generalized elementary functions}
\endsubhead
In this section we define the generalized elementary functions of one complex variable and the generalized elementary functions over a functional differential field. We also discuss a relation of these notions with generalized elementary extensions of differential fields. First we'll present needed lists of basic functions and of admissible operations.

\proclaim {List of basic elementary functions}

1. All complex constants and an independent variable $x$.

2. The exponential, the logarithm, and the power $x^\alpha$ where $\alpha$ is any constant.

3. The trigonometric functions sine, cosine, tangent, cotangent.

4. The inverse trigonometric functions arcsine, arccosine, arctangent, arccotangent.
\endproclaim
\proclaim{Lemma 2}
Basic elementary functions can be expressed through the exponentials and
the logarithms with the help of complex constants, arithmetic operations
and compositions.
\endproclaim
Lemma 2 can be considered as a  simple exercise. Its proof can be found in [2].

\proclaim {List of some classical  operations}

1. The operation of composition takes functions $f$,$g$ to the function $f\circ g$.

2.  The arithmetic operations take functions $f$, $g$ to the functions $f+g$, $f-g$, $fg$, and $f/g$.

3. The operation of differentiation takes function $f$ to the function $f'$.

4. The operation of solving algebraic equations takes functions $f_1,\dots,f_n$ to the function $y$ such that $y^n+f_1y^{n-1}+\dots+f_n=0$( the function $y$ is not quite uniquely determined by functions $f_1,\dots,f_n$ since an algebraic equation of degree $n$ can have $n$ solutions.
\endproclaim
\definition {Definition 4} The class of {\it generalized elementary functions of one variable} is defined by the following data:

List of basic functions: basic elementary functions.

List of admissible operations: Compositions, Arithmetic operations,  Differentiation, Operation of solving algebraic equations.
\enddefinition
\proclaim{Theorem 3}
A (possibly multivalued) function of one complex variable  belongs to the class of generalized elementary functions if and only if it belongs to some generalized elementary extension of the differential field of all rational functions of one variable.
\endproclaim

Theorem 3 follows from Lemma 2 (all needed  arguments can be found in [2]).

Let $K$ be  a functional differential field  consisting of meromorphic functional on a connected Riemann surface $U$ equipped with a meromorphic function $x$.

\definition {Definition 5} Class of {\it generalized elementary functions over a functional differential field $K$} is defined by the following data.

List of basic functions: all functions from the field $K$.

List of admissible operations: Operation of composition with a generalized elementary function $\phi$ that takes $f$ to $\phi\circ f$, Arithmetic operations,  Differentiation, Operation of solving algebraic equations.
\enddefinition

\proclaim{Theorem 4}
A (possibly multivalued) function on the Riemann surface $U$ belongs to the class of generalized elementary functions over a functional differential field $K$ if and only if it belongs to some generalize elementary extension of  $K$.
\endproclaim

Theorem 4 follows from Lemma 2 (all needed  arguments can be found in [2]).

\subhead{5. The Liouville's Theorem}
\endsubhead
In 1833 Joseph Liouville proved the following fundamental result.

\proclaim {Liouville's Theorem} {\it An  integral $y$ of a function $f$ from a functional differential field $K$ is a generalize4d elementary function over $K$
if and only if $y$ is representable in the form
$$
y(x)= \int\limits_{x_0}^x f(t)\,d t=r_0(x)+\sum\limits_{i=
1}^m\lambda_i\ln r_i(x), \tag 1
$$
where $r_0,\dots,r_m\in K$  and $\lambda_1,\dots,\lambda_m$ are complex constants.}
\endproclaim

For large classes of functions algorithms based on the Liouville's Theorem make it possible to either evaluate an integral or to prove that the integral cannot be ``evaluated in finite terms".

Let us  outline an inductive proof of the Liouville's Theorem.
\definition {Definition 6} A function $g$ is a {\it generalized elementary function of complexity $\leq k$} if there is a chain $K=F_0\subset F_1\subset\dots\subset F_k$ of functional differential fields such that $g\in F_k$ and for any $0\leq i<k$ either $F_{i+1}$ is a  finite  extension of  $F_{i}$, or $F_{i+1}$ is a pure transcendental extension of $F_i$ obtained by adjoining finitely many  exponentials, and  logarithms over $F_{i}$.
\enddefinition

We will prove the following induction hypothesis $I(m)$: {\it the Liouville's Theorem is true for every integral $y$ of complexity  $\leq m$ over any functional
differential field $K$}.
The statement $I(0)$ is obvious: if  $y\in K$, then $y=r_0\in K$. Now let  $y\,'\in K$ and $y\in F_k$. Since $y\,'\in F_1$,  by induction  $y=R_0 + \sum_{i=1}^q\lambda_i\ln R_i,$
where $R_0$, $R_1$, $\dots$, $R_q \in F_1$. We need to show that $y$ is representable in the form (1) with $r_0,\dots,r_m\in F_0=K$ We have the following two cases to consider:

1. $F_1$  is  a finite  extension of $F_0=K$. The statement of induction hypothesis in that case was proved by Abel and is called the Abel's Theorem. We will present its proof in the section 6.

2. $F_1$  a pure transcendental extension of $F_0=K$  obtained by adjoining exponentials and logarithms  over  $K$. We will deal with this case in section 7.

\head{6.  algebraic case}
\endhead
In the section 6.1 we  discuss finite extensions of differential fields. In the section 6.2 we  present a proof of the Abel's Theorem.

\subhead  {6.1.  An algebraic extension of a functional differential field}
\endsubhead
Let
$$
P(z)=z^n+a_1z^{n-1}+\dots+a_n \tag 2
$$
be an
irreducible polynomial over  $K$,
$P\in K[z]$.
Suppose that a functional differential field $F$ contains $K$ and a root $z$ of $P$.

\proclaim{Lemma 5}
 The field $K(z)$ is stable under the differentiation.

\endproclaim
\demo{Proof} Since  $P$ is irreducible over $K$, the polynomial
  $\frac{\partial P}{\partial z}$ has no common roots with $P$ and is different from zero
  in the field $K[z]/(P)$. Let $M$ be a polynomial  satisfying a congruence $M\frac{\partial P}{\partial z}\equiv -
  \frac{\partial P}{\partial x}\pmod P$.
  Differentiating the identity $P(z)=0$ in the field $F$,
we obtain that $\frac{\partial P}{\partial z}(z)z'+
\frac{\partial P}{\partial x}(z)=0$, which implies that $z'=M(z)$.
Thus the derivative of the element $z$ coincides with the value at $z$ of
a polynomial $M$. Lemma 5 follows from this fact.

\enddemo

Let $K\subset F$ and $\hat K\subset \hat F$ be  functional differential fields, and
$P$, $\hat P$
irreducible polynomials over $K$,
$\hat K$ correspondingly.
Suppose that  $F$, $\hat F$ contain  roots $z$, $\hat z$
of  $P$, $\hat P$.

\proclaim {Theorem 6} Assume that there is an isomorphism  $  \tau:K\rightarrow \hat K$ of differential fields $K$, $\hat K$ which maps coefficients of the polynomial $P$ to the corresponding coefficients of the polynomial $\hat P$. Then $\tau$ can be extended in a unique way to the differential isomorphism $ \rho:K(z)\rightarrow \hat K(\hat z)$.
\endproclaim

Proof of  Theorem 6 could be obtain by  the  arguments  used in the proof of Lemma~5.

\subhead {6.2. Induction hypothesis for an algebraic extension}
\endsubhead
Let $z_1,\dots,z_n$ be the roots of the polynomial  $P$ given by (2) and let $F_1=K\langle z_1 \rangle$.  Assume that there is an element $y_1\in F_1$, such that  $y_1'\in K$, $M_i\in K[x]$ and $y_1'$ is representable in the form
$$
y_1'=\sum_{i=1}^q\lambda_i\frac{(M_i(z_1))'}{M_i(z_1)}+(M_0(z_1))'. \tag 3
$$

\proclaim {Abel's Theorem} Under the above assumptions the element $y_1'$ is representable in the form (1) with polynomials  $M_i$ independent of $z_1$, i.e. with $M_0,M_1,\dots,M_q\in K$.
\endproclaim

\demo {Proof} Let $y_1$ be equal to $Q(z_1)$ where $Q\in K[z]$. For any $1\leq j\leq n$ let $y_j$ be the element $Q(z_j)$. According to Theorem 6 the identity (3) implies the identity

$$
y_j'=\sum_{i=1}^q\lambda_i\frac{(M_i(z_j))'}{M_j(z_1)}+(M_0(z_j))'. \tag 4
$$

Since $y'_1\in K$ we obtain $n$ equalities  $y'_1=\dots=y'_n$. To complete the proof it is enough to
take the arithmetic mean of  $n$ equalities (4). Indeed  the elements
$\tilde M_i=\prod_{1\leq k\leq n}M_i(z_k)$ and $\tilde M_0=\sum_{1\leq k\leq n}M_0(z_k)$ are symmetric functions in the roots of the polynomial $P$ thus $\tilde M_0, \dots \tilde M_q\in K$.
\enddemo

\remark{Remark}
The proof uses implicitly the Galois group $G$ of the  splitting field of the polynomial $P$ over the field $K$. The group $G$ permutes the roots $y_1,\dots,y_n$ of $P$. The  element $\tilde M_i=\prod_{1\leq k\leq n}M_i(z_k)$ and $\tilde M_0=\sum_{1\leq k\leq n}M_0(z_k)$ are invariant under the action of $G$  thus they belong to the field $K$.
\endremark

\head{III.  pure transcendental  case}
\endhead

In this chapter we prove induction hypothesis in the pure transcendental case. First we will state the corresponding  Theorem 7 and will outline its proof.

Let $F_1$ be a functional differential field obtained by  extension of the functional differential field $K$   by adjoining algebraically independent over $K$ functions
$$
y_1=\ln a_1, \dots, y_k=\ln a_k, z_1=\exp b_1,\dots,z_n=\exp b_n \tag 5
$$
where $a_1,$ $\dots,$ $a_k,$ $b _1,$ $\dots,$ $b_k$ are some functions from $K$. We will assume that $F_1$ consists of meromorphic functions on a connected Riemann surface $U$ and the differentiation in $K_1$ using a meromorphic function $x$ on $U$. Let $X$ be the manifold $U\times G$ where $G=\Bbb C^k\times (\Bbb C^*)^n$. Consider a map $\gamma:U\rightarrow \Bbb C^k\times (\Bbb C^*)^n$given by formula
$$
\gamma(p)=y_1(p),\dots, y_k(p),z_1(p)\dots,z_n(p)
$$
where the functions $y_i$, $z_j$ are defined by (5).

Let $X$ be the product $U\times (\Bbb C)^k\times (\Bbb C^*)^n$. Denote by $\Gamma\subset X$ the graph of the map~$\gamma$. Consider a germ $\Phi$ of a complex valued function at the point $a\in X$.

\definition {Definition 7}
We say that  $\Phi$ is a {\it logarithmic type germ} if $\Phi$ is representable in the form $\Phi_a=R_0+\sum_{i=1}^q\lambda_i\ln R_i,$
where $R_i$ are germs at the point $a\in X$ of rational functions of $(y_1,\dots,y_k, z_1,\dots,z_n)$ with coefficients in $K$ and $\lambda_j$ are complex numbers.
\enddefinition

\proclaim { Theorem 7}  Let $\Phi$ be a logarithmic type germ at a point $a=(p_0,\gamma (p_0))\in \Gamma$. Then the germ of the function $\Phi (p,\gamma(p))$ at the point $p_0\in U$  is a germ of an  integral over $K$ if and only if $\Phi$ is representable in the following form
$$\Phi(p,y,z)=\Phi(p, \gamma(p_0))+\sum_1^k c_i (y_i -y_i(p_0))+\sum_1^nt_j\ln \frac{z_j}{z_j(p_0)}\tag 6$$ where $r_0$ is a germ of a function from the field $K$ and $c_i,t_j$ are complex constants.
\endproclaim

Theorem 7 proves induction hypothesis in the pure transcendental case. Indeed the germ $\Phi(p,\gamma(p_0))$ given by (6) is a germ of a function from the field $K$ and according to (5) the  identities $c_iy_i=c_i\ln a_i$, $t_j\ln z_j=t_jb_j$ hold. We split the claim of Theorem 7 into two parts.

First we consider the natural action of the group $G=(\Bbb C^k)\times (\Bbb C^*)^n$ on $X=U\times G$ and we  describe all germs of closed 1-forms locally invariant under this action.  Corollary 11 claims that each such  1-form is a differential of a function representable in the form (6).

Second we show that if the germ $\Phi$ satisfies the conditions of Therem then the germ $d\Phi$ is locally invariant under the action of the group $G$ (see Theorem 16).

\subhead {7.1. Locally invariant closed 1-forms}
\endsubhead
Let $G$ be a connected Lie group acting by diffeomorphisms on a manifold $X$. Let $\pi:G\rightarrow Diff (X)$ be a corresponding homomorphism from $G$ to the group $Diff (X)$ of diffeomorphisms of $X$. For a vector $\xi$ from the Lie algebra $\Cal G$ of $G$ the action $\pi$ associates the vector field $V_\xi$ on $X$. The germ $\omega_{x_0}$ at a point $x_0\in X$ of a differential form $\omega$ on $X$ is {\it locally invariant under the action $\pi$} if for any $\xi\in \Cal G$ the Lie derivative $L_{V_\xi}\omega$ is equal to zero.

\proclaim{Lemma 8} The germ of the differential $d \varphi_{x_0}=\omega_{x_0}$ of a smooth function $\varphi$  is locally invariant under the action $\pi$ if and only if for each  $\xi\in \Cal G$ the Lie derivative $L_{V_\xi}\varphi$ is a constant $M(\xi)$ (which depends on $\xi$).
\endproclaim
\demo{Proof}  Applying ``Cartan's magic formula"
$
L_{V_\xi}\omega =i_{X}d\omega +d(i_{X}\omega )
$
we obtain that $L_{V_\xi}\omega=0$ if and only if $d(L_{V_\xi}\varphi)=0 $ which means that  $L_{V_\xi}\varphi$ is constant.
\enddemo

The following Theorem characterizes locally invariant  closed 1-forms  more explicitly.

\proclaim{Theorem 9}  The germ of the differential $d \varphi_{x_0}=\omega_{x_0}$ of a smooth complex valued function $\varphi$  is locally invariant under the action $\pi$ if and only if there exists a local homomorphism $\rho$ of $G$ to the additive group  $\Bbb C$ of complex numbers such that for any $g\in G$ in a neighborhood of the identity  the following relation  holds:
$$
\varphi(\pi(g)x_0)=\varphi(x_0)+\rho(g).
$$

\endproclaim

\demo{Proof} For $\xi\in \Cal G$  the Lie derivative $L_{V_\xi}\varphi$ is constant $M(\xi)$  by Lemma 8. Let us show that for $\xi\in [\Cal G ]$ where $[\Cal G ]$ is the commutator  of $\Cal G$ the constant $M(\xi)$ equals to zero. Indeed if $\xi=[\tau,\rho]$ then
$$
L_{V_\xi} \varphi= L_{V_\rho}L_{V_\tau}\varphi-L_{V_\tau}L_{V_\rho}\varphi=L_{V_\rho}M(\tau)-L_{V_\tau}M(\rho)=0.
$$
Thus the linear function $M:\Cal G \rightarrow \Bbb C$  mapping $\xi$ to $M(\xi)$ provides a homomorphism of $\Cal G $ to the Lie algebra of the additive group $\Bbb C$ of complex numbers. Let $\rho$ be the local homomorphism of $G$ to $\Bbb C$ corresponding to the homomorphism $M$.

Consider a function $\phi$  on a neighborhood of the identity  in $G$ defined by the following formula: $\phi(g)=\varphi(x_0)+\rho(g)$. By definition on a neighborhood of identity the function $\phi$ has the same differential as the function $\varphi(\pi(g)x)$. Values  of these functions at the identity  are equal to $\varphi(x_0)$. Thus these functions are equal.
\enddemo

Assume that  $X=U\times G$ where $U$ is a  manifold  and an action $\pi$ is given by the formula $\pi(g) (x,g_1)= (x, gg_1)$. Applying Theorem 9 to this action  we obtain the following corollary.

\proclaim {Corollary 10} If  germ of  differential $d \varphi=\omega$ of a smooth complex valued function $\varphi$ at a point $(x_0,g_0)\in U\times G$ is locally invariant under the action $\pi$ then  in a neighborhood of the point $(x_0,g_0)$ the following identity holds:
$$
\varphi(x, g)=\varphi(x, g_0)+\rho(gg_0^{-1}).\tag 7
$$
where $\rho$ is a local homomorphism of $G$ to the additive group of complex numbers.

\endproclaim

\demo{Proof} Follows from Theorem 9 since the element $gg_0^{-1}$ maps the point $(x,g_0)$ to the point $(x,g)$.
\enddemo

Let $G$ be the group $\Bbb C^k\times (\Bbb C^*)^n$ where $\Bbb C$ and $\Bbb C^*$ are additive and mulplicative group of complex numbers. We will consider the group $\Bbb C^k\times (\Bbb C^*)^n$ with coordinate functions   $(y,z)=(y_1,\dots, y_k,z_1,\dots, z_n)$ assuming that $z_1\cdot\dots\cdot z_n\neq 0$.
\proclaim {Corollary 11} If in the assumptions of Corollary 1 for $G=\Bbb C^k\times (\Bbb C^*)^n$ in a neighborhood of $(x_0,y_0,z_0)\in U\times(\Bbb C^k\times (\Bbb C^*)^n )$ the following identity holds
$$
\varphi(x,y,z)=\varphi(x,y_0,z_0)+\sum_{1\leq i\leq k} \lambda_i(y_i-(y_0)_{i}) +\sum_{1\leq j \leq n} \mu_j\ln \frac{z_j}{(z_0)_{j}}
$$
where $\lambda_1,\dots,\lambda_k, \mu_1\dots,\mu_n$ are complex constants.

\endproclaim
\demo{Proof} Follows from (7) since any local homomorphism $\rho$ from the group $\Bbb C^k\times (\Bbb C^*)^n$ to the additive group of complex numbers can be given by formula
 $$
 \rho(y_1,\dots,y_k,z_1,\dots,z_n)=\sum _{1\leq i \leq k}\lambda_1 y_i+\sum_{1\leq j \leq n}\mu_j\ln z_j
 $$
 where $\lambda_i$ and $\mu_j$ are complex constants.
\enddemo

\subhead {7.2. Vector field associated to a logarithmic-exponential extension}
\endsubhead
We use the notations introduced  in  the section 7. Let $G$ be the group $\Bbb C^k\times (\Bbb C^*)^n$ and let $X$ be the product $U\times G$ consider the map $\gamma:U\rightarrow  \Bbb C^k\times (\Bbb C^*)^n$ given by the following formula:

$$y_1=\ln a_1,\dots,y_k=\ln a_k, z_1=\exp b_1,\dots,z_n=\exp b_n.\tag 8$$

The map $\gamma$ satisfies  the following differential relation:
$$
  d\gamma= da_1/a_1,\dots,da_k/a_k,z_1d b_1,\dots,z_ndb_n
$$
\definition {Definition 8} Let $V$ be a meromorphic vector field on $X$ defined by the following conditions. If $V_a$ is the value of $V$ at the point $a=(p,y_1,\dots, y_k,z_1,\dots,z_n)\in X$ then $\langle dx, V_a\rangle =1$, $\langle dy_i, V_a\rangle =a_i'/a_i(p)$ for $1\leq i\leq k$, $\langle dz_j, V_a\rangle =b_j'/b_i(p)$ for $1\leq j\leq n$
\enddefinition

Vector field $V$ is regular on $U^0\times G$ where $U^0$ is an open subset in $U$ which does not contain zeros and poles of the functions $a_1,\dots, a_k$ and poles of functions $b_1,\dots,b_n$ poles and zeros and poles of the 1-form $dx$. By construction the graph $\Gamma=(p,\gamma(p)\subset X$ of the map $\gamma$ is an integral curve for differential equation on $X$ defined by the vector field $v$.

The following lemmas are obvious.

\proclaim {Lemma 12}  The vector field $V$  is invariant under the action $\pi$ on $X$. For each element $g\in G$ the curve  $g \Gamma\subset X$ of the graph $\Gamma$ of $\gamma$ is an integral curve for $V$.
\endproclaim

\proclaim {Lemma 13}  The field $K(y,z)$ of rational functions in  $y_1,\dots,y_k,z_1,\dots, z_n $ over the field $K$ is invariant under the action $\pi$ on $X$. For each vector $\xi\in \Cal G$ in the Lie algebra $\Cal G$ of $G$ the Lie derivative $L_{V_\xi} R$ of $R\in K(y,z)$ belongs to $K(y,z)$.
\endproclaim

\subhead {7.3. Pure transcendental logarithmic exponential extension}
\endsubhead
We will assume below that the components (8) of $\gamma$ are algebraically independent over $K$.

{\bf Lioville's principal.} If a polynomial $P\in K[y_1,\dots,y_k,z_1,\dots,z_n]$ vanishes on the graph $\Gamma\subset X$ of the map $\gamma$ then $P$ is identically equal to zero.

\demo {Proof} If $P$ is not identically equal to zero then the components of $\gamma$ are algebraically dependent over the field $K$.
\enddemo

\proclaim{Theorem 14} The extension $K\subset F_1$ is isomorphic to the extension of $K$ by the field of rational functions $K(y,z)$ in $(y_1,\dots,y_k,z_1,\dots,z_n)$ over $K$ considered as the field of functions on $X$  equipped with the differentiation sending $f\in K(y,z)$ to the Lie derivative $L_Vf$ with respect to the vector field $V$ introduced in definition 8.
\endproclaim
\demo{Proof} By assumption components (8) of the map $\gamma$ are algebraically independent over $K$ thus each function from the extension obtained by adjoining to  $K$ by these components is representable in the unique way as a rational function from $K(y,z)$. By definition the derivatives of the components (8) are coincide with their Lie derivatives  with respect to the vector field $V$.
\enddemo

The action $\pi$ of the group $G=\Bbb C^k\times(\Bbb C^*)^n$ on $X$ induces the action $\pi^*$ of $G$ on the space of functions on $X$ containing the field  $K(y,z)$. The vector field $V$ is invariant under the action $\pi$. Thus $\pi^*$ acts on $K(y,z)\sim F_1$ by  differential automorphisms. Easy to see that a function $f\in K(y,z)$ is fixed under the action $\pi^*$   if and only if  $f\in K$, i.e. the group $G$ is isomorphic to the {\it differential Galois group} of the extension $K\subset F_1$. We proved the following result

\proclaim{Theorem 15} The differential Galois group of the extension $K\subset F_1$ is isomorphic to the group $G$. The Galois group  is induced on the differential field $K(y,z)$ with the differentiation given by Lie derivative with respect to the field $V$   by the action of $G$ on the manifold $X=U\times \Bbb C^k\times (\Bbb C^*)^n$.
\endproclaim

Now we are ready to complete inductive proof of the Liouvvile's Theorem.

\proclaim {Theorem 16}  Let $\Phi$ be a logarithmic type germ at a point $a=(p_0,\gamma (p_0))\in \Gamma\subset X$. If the germ of the function $\Phi (p,\gamma(p))$ on $U$ at the point $p_0\in U$  is a germ of an  integral $f$ over $K$ then the germ of the differential $d\Phi$ at the point $a\in X$ is locally invariant under the action $\pi$ on $X$.

\endproclaim

\demo {Proof} By the assumption of Theorem the restriction of the function $(L_V\Phi-f)$  on $\Gamma$  is equal to zero. Since the function $(L_V\Phi-f)$ belongs to the field $K(y,z)$  the function $(L_V\Phi-f)$ by the Liouville's principal is equal to zero identically on $X$. In particular it is equal to zero  on the integral curve $g \Gamma$ the vector field  $V$, where $g$ is an element of the group $G$. Thus  the restrictions of function  $L_V\pi(g)^*(\Phi-f)$ to $\Gamma$ equals to zero. Since the function $f$ is invariant under the action $\pi^*$ we obtain that the restriction on $\Gamma$ of $L_V(\Phi-\pi^*(g)\Phi)$ is equal to zero. Differentiating this identity we obtain that for any $\xi\in \Cal G$ the restriction on $\Gamma$ of $L_V(L_{V_\xi}\Phi$ equals to zero.
Thus on $\Gamma$ the function $L_\xi\Phi$ is constant.  Lemma 13 implies  that the function $L_{V_\xi}$ belongs to the field $K(x,y)$. Thus  the function$L_{V_\xi}$ is a constant on $X$ by Liouville's principal. Thus the 1-form $d\Phi$ is locally invariant under the action $\pi$  by Lemma 8. Theorem 16 is proved.
\enddemo

Thus we complete proof of Theorem 7 and  the inductive proof of Liouville's Theorem.

\end

\tc  Abstrakt. Soglasno teoreme Liuvillya neopredelenny\u i integral e1lementarno\u i funktsii kak pravilo ne yavlyaet\-sya elementarno\u i funktsie\u i. V nastoyashche\u i statp1e my obsuzhdaem tochnuyu formulirovku i dokazatelp1stvo ep1to\u i teoremy. Differentsialp1naya gruppa Galua  rassshireniya, poluchennogo prisoedineniem integrala ne pozvolyaet opredelitp1 beret\-sya ep1tot integral v e1lementarnykh funktsiyakh ili net.  Tem ne menee teoremu Liuvillya mozhno dokazatp1 ispolp1zuya differentsialp1nye gruppy Galua. Pervy\u i shag v e1tom napravlenii byl sdelan Abelem. E1tot shag svyazan s konechnymi rasshireniyami i ikh konechnymi gruppami Galua. Znachitelp1naya chastp1 statp1i posvyashchena sleduyushchemu shagu, svyazannomu s chisto transtcendentnymi rasshireniyami i ikh differentsialp1nymi gruppami Galua, kotorye predstavlyayut sobo\u i svyaznye algebraicheskie gruppy Li. Ideya dokazatelp1stva voskhodit k Zh. Liuvillyu i  Zh. F. Rittu.
\end